\documentclass[11pt]{article} 
\usepackage{amsfonts,amsmath,latexsym,amssymb,mathrsfs,amsthm,comment}
\usepackage{slashbox}
\usepackage{caption}

\let\OLDthebibliography\thebibliography
\renewcommand\thebibliography[1]{
  \OLDthebibliography{#1}
  \setlength{\parskip}{1pt}
  \setlength{\itemsep}{1pt plus 0.1ex}
}

\def\numberlikeadb{\global\def\theequation{\thesection.\arabic{equation}}}
\numberlikeadb
\newtheorem{theorem}{Theorem}[section]

\newtheorem{corollary}[theorem]{Corollary}

\newtheorem{remark}[theorem]{Remark}

\usepackage{color}

\usepackage{lscape}
\usepackage{caption}
\usepackage{multirow}
\allowdisplaybreaks
\begin{document}

\title{On the cumulative distribution function of the variance-gamma distribution
}
\author{Robert E. Gaunt\footnote{Department of Mathematics, The University of Manchester, Oxford Road, Manchester M13 9PL, UK}  
}

\date{} 
\maketitle

\vspace{-10mm}

\begin{abstract}We obtain exact formulas for the cumulative distribution function of the variance-gamma distribution, as infinite series involving the modified Bessel function of the second kind and the modified Lommel function of the first kind. From these formulas, we deduce exact formulas for the cumulative distribution function of the product of two correlated zero mean normal random variables.

\end{abstract}


\noindent{{{\bf{AMS 2010 Subject Classification:}}} Primary 60E05; 62E15

\section{Introduction}

The variance-gamma (VG) distribution with parameters $\nu > -1/2$, $0\leq|\beta|<\alpha$, $\mu \in \mathbb{R}$, denoted by $\mathrm{VG}(\nu,\alpha,\beta,\mu)$, has probability density function (PDF)
\begin{equation}\label{vgpdf} p(x) = M \mathrm{e}^{\beta (x-\mu)}|x-\mu|^{\nu}K_{\nu}(\alpha|x-\mu|), \quad x\in \mathbb{R},
\end{equation}
where the  normalising constant is given by
\[M=M_{\nu,\alpha,\beta}=\frac{(\alpha^2-\beta^2)^{\nu+1/2}}{\sqrt{\pi}(2\alpha)^\nu \Gamma(\nu+1/2)},\]
and $K_\nu(x)$ is a modified Bessel function of the second kind (see Appendix \ref{appa} for a definition). The parameters have the following interpretation: $\nu$ is a shape parameter, $\alpha$ is a scale parameter, $\beta$ is a skewness parameter, and $\mu$ is a location parameter. Other names include the Bessel function distribution \cite{m32}, the McKay Type II distribution \cite{ha04} and the generalized Laplace distribution \cite[Section 4.1]{kkp01}. Alternative parametrisations are given in \cite{gaunt vg,kkp01,mcc98}. Interest in the VG distribution dates as far back as 1929 in which the VG PDF (\ref{vgpdf}) arose as the
PDF of the sample covariance for a random sample drawn from a bivariate normal population \cite{p29}. The VG distribution was introduced into the financial literature in the seminal works \cite{mcc98,madan}, and has recently found application in probability theory as a natural limit distribution \cite{aet21,gaunt vg}.
Further application areas and distributional properties can be found in the survey \cite{vgsurvey} and the book \cite{kkp01}.

In this paper, we fill in an obvious gap in the literature by deriving exact formulas for the cumulative distribution function (CDF) of the VG distribution that hold for the full range of parameter values. Our formulas are expressed as infinite series involving the modified Bessel function of the second kind and the modified Lommel function of the first kind (defined in Appendix \ref{appa}). 
Despite being widely used in financial modelling and other applications areas, exact formulas had only previously been given for the symmetric case $\beta=0$ \cite{jp} and for the case $\nu\in\{1/2,3/2,5/2,\ldots\}$ \cite{nsg07}, in which case the modified Bessel function $K_\nu(x)$ in the PDF (\ref{vgpdf}) takes an elementary form; see equation (\ref{special}). 

As the product of two correlated zero mean normal random variables, and more generally the sum of $n\geq1$ independent copies of such random variables are VG distributed \cite{gaunt prod}, we immediately deduce exact formulas for the CDFs of these distributions. These distributions also have numerous applications, dating back to the work of \cite{craig} in 1936; for an overview of application areas and distributional properties see \cite{gaunt22}. Since the work of \cite{craig}, the problem of finding the exact PDF of these distributions has received much interest; see \cite{np16} for an overview of the contributions in the literature. We thus contribute to the next natural problem of finding exact formulas for the CDF. Formulas for the CDF for the case $n\geq2$ is an even integer have been obtained by \cite{gaunt22} (in this case the PDF takes an elementary form, which is again a consequence of equation (\ref{special})). In this paper, we obtain formulas for the CDF that hold for all $n\geq1$, which includes the important case $n=1$ for the distribution of a single product of two correlated zero mean normal random variables. 


\section{Results and proofs}

The following theorem is the main result of this paper. Let $F_X(x)=\mathbb{P}(X\leq x)$ denote the CDF of $X\sim \mathrm{VG}(\nu,\alpha,\beta,\mu)$. Also, for $\mu\geq\nu>-1/2$, let
\begin{align}\label{gmn}G_{\mu,\nu}(x)&=x\big(K_{\nu}( x)\tilde{t}_{\mu-1,\nu-1}( x)+K_{\nu-1}( x)\tilde{t}_{\mu,\nu}( x)\big),\\
\label{gmn2}\tilde{G}_{\mu,\nu}(x)&=1-G_{\mu,\nu}(x),
\end{align}
where $\tilde{t}_{\mu,\nu}(x)$ is a normalisation of the modified Lommel function of the first kind $t_{\mu,\nu}(x)$, defined in Appendix \ref{appa}. In interpreting the formulas in the theorem, it should be noted that, for fixed $\mu\geq\nu>-1/2$, $G_{\mu,\nu}(x)$ ($\tilde{G}_{\mu,\nu}(x)$) is an increasing (decreasing) function of $x$ on $(0,\infty)$ satisfying $0<G_{\mu,\nu}(x)<1$ and $0<\tilde{G}_{\mu,x}(x)<1$ for $x>0$ (see Appendix \ref{appa}). One of the formulas in the theorem is also expressed in terms of the hypergeometric function, which is defined in Appendix \ref{appa}. We also let $\mathrm{sgn}(x)$ denote the sign function, $\mathrm{sgn}(x)=1$ for $x>0$, $\mathrm{sgn}(0)=0$, $\mathrm{sgn}(x)=-1$ for $x<0$. 

\begin{theorem}\label{thm1} Let $X\sim \mathrm{VG}(\nu,\alpha,\beta,\mu)$, where $\nu > -1/2$, $0\leq|\beta|<\alpha$, $\mu\in\mathbb{R}$. Then, for $x\geq\mu$,
\begin{align}\label{eq1}F_X(x)=1&-\frac{(1-\beta^2/\alpha^2)^{\nu+1/2}}{2\sqrt{\pi}\Gamma(\nu+1/2)}\sum_{k=0}^\infty\frac{1}{k!}\bigg(\frac{2\beta}{\alpha}\bigg)^k\Gamma\bigg(\frac{k+1}{2}\bigg)\Gamma\bigg(\nu+\frac{k+1}{2}\bigg)\tilde{G}_{\nu+k,\nu}(\alpha (x-\mu)),
\end{align}
and, for $x<\mu$,
\begin{align}\label{eq2}F_X(x)&=\frac{(1-\beta^2/\alpha^2)^{\nu+1/2}}{2\sqrt{\pi}\Gamma(\nu+1/2)}\sum_{k=0}^\infty\frac{(-1)^k}{k!}\bigg(\frac{2\beta}{\alpha}\bigg)^k\Gamma\bigg(\frac{k+1}{2}\bigg)\Gamma\bigg(\nu+\frac{k+1}{2}\bigg)\tilde{G}_{\nu+k,\nu}(-\alpha (x-\mu)).
\end{align}
Moreover, the following formula is valid for all $x\in\mathbb{R}$:
\begin{align}F_X(x)&=\frac{1}{2}-\frac{\Gamma(\nu+1)}{\sqrt{\pi}\Gamma(\nu+1/2)}\frac{\beta}{\alpha}\bigg(1-\frac{\beta^2}{\alpha^2}\bigg)^{\nu+1/2}{}_2F_1\bigg(1,\nu+1;\frac{3}{2};\frac{\beta^2}{\alpha^2}\bigg)\nonumber\\
\label{eq3}&\quad+\frac{(1-\beta^2/\alpha^2)^{\nu+1/2}}{2\sqrt{\pi}\Gamma(\nu+1/2)}\sum_{k=0}^\infty\frac{(\mathrm{sgn}(x))^{k+1}}{k!}\bigg(\frac{2\beta}{\alpha}\bigg)^k\Gamma\bigg(\frac{k+1}{2}\bigg)\Gamma\bigg(\nu+\frac{k+1}{2}\bigg)\times\nonumber\\
&\quad\times G_{\nu+k,\nu}(\alpha |x-\mu|).
\end{align}
\end{theorem}

\begin{remark}1. Let $X\sim \mathrm{VG}(\nu,\alpha,\beta,\mu)$, where $\nu > -1/2$, $0\leq|\beta|<\alpha$, $\mu\in\mathbb{R}$. The probability $\mathbb{P}(X\leq\mu)$ takes a particularly simple form:
\begin{equation*}\mathbb{P}(X\leq\mu)=\frac{1}{2}-\frac{\Gamma(\nu+1)}{\sqrt{\pi}\Gamma(\nu+1/2)}\frac{\beta}{\alpha}\bigg(1-\frac{\beta^2}{\alpha^2}\bigg)^{\nu+1/2}{}_2F_1\bigg(1,\nu+1;\frac{3}{2};\frac{\beta^2}{\alpha^2}\bigg).
\end{equation*}
We used \emph{Mathematica} to calculate this probability for the case $\alpha=1$ and $\mu=0$, for a range of values of $\nu$ and $\beta$; the results are reported in Table \ref{table1}. We only considered positive values of $\beta$ due to the fact that if $Y\sim \mathrm{VG}(\nu,1,\beta,0)$, then $-Y\sim\mathrm{VG}(\nu,1,-\beta,0)$ (see \cite[Section 2.1]{vgsurvey}). We observe from Table \ref{table1} that the probability $\mathbb{P}(Y\leq0)$ decreases as the skewness parameter $\beta$ increases and as the shape parameter $\nu$ increases. 

\vspace{1mm}

\noindent 2. The CDF takes a simpler form when $\beta=0$. Suppose that $X\sim\mathrm{VG}(\nu,\alpha,0,\mu)$. Then applying (\ref{struve}) to (\ref{eq3}) yields the following formula for the CDF of $X$: for $x\in\mathbb{R}$,
\begin{align*}F_X(x)=\frac{1}{2}+\frac{\alpha(x-\mu)}{2}\bigg[K_{\nu}(\alpha|x-\mu|)\mathbf{L}_{\nu-1}(\alpha|x-\mu|)+\mathbf{L}_{\nu}(\alpha|x-\mu|)K_{\nu-1}(\alpha|x-\mu|)\bigg],
\end{align*}
where $\mathbf{L}_\nu(x)$ is a modified Struve function of the first kind (see \cite[Chapter 11]{olver} for a definition and properties). Other formulas for the special case $\beta=0$ are given by \cite{jp}.
\end{remark}

\begin{table}[h]
\centering
\caption{\footnotesize{$\mathbb{P}(Y\leq0)$ for $Y\sim\mathrm{VG}(\nu,1,\beta,0)$.}}
\label{table1}
{\normalsize
\begin{tabular}{|c|rrrrrrr|}
\hline
 \backslashbox{$\beta$}{$\nu$}       &    $-0.25$ &    0 &    0.5 & 1 &  2 &  3 &    5   \\
 \hline
 0.05 & 0.4905 & 0.4841 & 0.4750 & 0.4682 & 0.4576 & 0.4492 & 0.4356 \\
 0.1  & 0.4809 & 0.4681 & 0.4500 & 0.4364 & 0.4155 & 0.3990 & 0.3726 \\
0.25  & 0.4516 & 0.4196 & 0.3750 & 0.3425 & 0.2944 & 0.2582 & 0.2050 \\
0.5  & 0.3978 & 0.3333 & 0.2500 & 0.1955 & 0.1266 & 0.0852  & 0.0409  \\
0.75  & 0.3271 & 0.2301 & 0.1250 & 0.0721 & 0.0261 & 0.0100 & 0.0016  \\
  \hline
\end{tabular}}
\end{table}

\begin{proof}To ease notation, we set $\mu=0$; the general case follows from a simple translation. Suppose first that $x\geq0$. Using the formula (\ref{vgpdf}) for the VG PDF, the power series expansion of the exponential function, and interchanging the order of integration and summation gives that
\begin{align*}F_X(x)=1-M\int_x^\infty\mathrm{e}^{\beta t}t^\nu K_{\nu}(\alpha t)\,\mathrm{d}t=1-M\sum_{k=0}^\infty\frac{\beta^k}{k!}\int_x^\infty t^{\nu+k}K_\nu(\alpha t)\,\mathrm{d}t.
\end{align*}
Evaluating the integrals using the integral formula (\ref{kint}) 
yields the formula (\ref{eq1}).

Now suppose $x<0$. Arguing as before, we obtain that
\begin{align}F_X(x)&=M\int_{-\infty}^x \mathrm{e}^{\beta t}(-t)^\nu K_{\nu}(-\alpha t)\,\mathrm{d}t=M\sum_{k=0}^\infty\frac{\beta^k}{k!}\int_{-\infty}^x (-1)^k(-t)^{\nu+k}K_\nu(-\alpha t)\,\mathrm{d}t\nonumber\\
\label{hjk}&=M\sum_{k=0}^\infty\frac{(-\beta)^k}{k!}\int_{-x}^\infty y^{\nu+k}K_\nu(\alpha y)\,\mathrm{d}y,
\end{align}
and evaluating the integrals in (\ref{hjk}) using (\ref{kint}) yields the formula (\ref{eq2}).

We now derive formula (\ref{eq3}). Let $x\in\mathbb{R}$. Proceeding as before, we obtain that
\begin{align}F_X(x)&=F_X(0)+M\mathrm{sgn}(x)\int_0^x\mathrm{e}^{\beta t}|t|^\nu K_{\nu}(\alpha |t|)\,\mathrm{d}t\nonumber\\
&=F_X(0)+M\mathrm{sgn}(x)\sum_{k=0}^\infty\frac{\beta^k}{k!}\int_0^x (-1)^k|t|^{\nu+k}K_\nu(\alpha |t|)\,\mathrm{d}t\nonumber\\
\label{july}&=F_X(0)+M\sum_{k=0}^\infty\frac{\beta^k}{k!}(\mathrm{sgn}(x))^{k+1}\int_0^{|x|} t^{\nu+k}K_\nu(\alpha t)\,\mathrm{d}t.
\end{align}
The integrals in (\ref{july}) can be evaluated using the integral formula (\ref{kint2}), and it now remains to compute $F_X(0)$. 

Applying formula (\ref{eq2}) with $x=0$ and using that $\lim_{x\rightarrow0}\tilde{G}_{\nu+k,\nu}(x)=1$ (this is readily obtained by applying the limiting forms (\ref{Ktend0}) and (\ref{tend0})) yields
\begin{align}F_X(0)&=\frac{(1-\beta^2/\alpha^2)^{\nu+1/2}}{2\sqrt{\pi}\Gamma(\nu+1/2)}\sum_{k=0}^\infty\frac{(-1)^k}{k!}\bigg(\frac{2\beta}{\alpha}\bigg)^k\Gamma\bigg(\frac{k+1}{2}\bigg)\Gamma\bigg(\nu+\frac{k+1}{2}\bigg)\nonumber\\
\label{s3}&=\frac{(1-\beta^2/\alpha^2)^{\nu+1/2}}{2\sqrt{\pi}\Gamma(\nu+1/2)}\big(S_1+S_2),
\end{align}
where
\begin{align*}S_1&=\sum_{k=0}^\infty\frac{1}{(2k)!}\bigg(\frac{2\beta}{\alpha}\bigg)^{2k}\Gamma\bigg(k+\frac{1}{2}\bigg)\Gamma\bigg(\nu+k+\frac{1}{2}\bigg),\\
S_2&=-\sum_{k=0}^\infty\frac{1}{(2k+1)!}\bigg(\frac{2\beta}{\alpha}\bigg)^{2k+1}k!\Gamma(\nu+k+1).
\end{align*}
On calculating $(2k)!=\Gamma(2k+1)$ using the formula $\Gamma(2x)=\pi^{-1/2}2^{2x-1}\Gamma(x)\Gamma(x+1/2)$ (see \cite[Section 5.5(iii)]{olver}) and then applying the standard formula $(u)_k=\Gamma(u+k)/\Gamma(u)$, we obtain
\begin{align}\label{s1}S_1=\sqrt{\pi}\Gamma(\nu+1/2)\sum_{k=0}^\infty\frac{(\nu+1/2)_k}{k!}\bigg(\frac{\beta}{\alpha}\bigg)^{2k}=\frac{\sqrt{\pi}\Gamma(\nu+1/2)}{(1-\beta^2/\alpha^2)^{\nu+1/2}},
\end{align}
where we evaluated the sum using the generalized binomial theorem. Using similar considerations, we can express $S_2$  in the hypergeometric form (\ref{gauss}), which yields
\begin{align}\label{s2}S_2=-\frac{2\beta}{\alpha}\Gamma(\nu+1){}_2F_1\bigg(1,\nu+1;\frac{3}{2};\frac{\beta^2}{\alpha^2}\bigg).
\end{align}
Plugging formulas (\ref{s1}) and (\ref{s2}) into (\ref{s3}) now yields formula (\ref{eq3}).
\end{proof}

Now, we let $(U,V)$ be a bivariate normal random vector having  zero mean vector, variances $(\sigma_U^2,\sigma_V^2)$ and correlation coefficient $\rho$. Let $Z=UV$ be the product of these correlated normal random variables, and let $s=\sigma_U\sigma_V$.  We also introduce the mean $\overline{Z}_n=n^{-1}(Z_1+Z_2+\cdots+Z_n)$, where $Z_1,Z_2,\ldots,Z_n$ are independent copies of $Z$. It was noted by \cite{gaunt thesis} that $Z$ is VG distributed,  and more generally it was shown by \cite{gaunt prod} that
\begin{equation}\label{vgrep}\overline{Z}_n\sim\mathrm{VG}\bigg(\frac{n-1}{2},\frac{n}{s(1-\rho^2)},\frac{n\rho}{s(1-\rho^2)},0\bigg).
\end{equation}
On combining (\ref{vgrep}) with (\ref{eq1}), (\ref{eq2}) and (\ref{eq3}), we obtain the following formulas for the CDF of $\overline{Z}_n$; formulas for the CDF of $Z$ are obtained by letting $n=1$. 

\begin{corollary}Let the previous notations prevail. Then, for $x\geq0$,
\begin{align*}F_{\overline{Z}_n}(x)&=1-\frac{(1-\rho^2)^{n/2}}{2\sqrt{\pi}\Gamma(n/2)}\sum_{k=0}^\infty\frac{(2\rho)^k}{k!}\Gamma\bigg(\frac{k+1}{2}\bigg)\bigg(\frac{n+k}{2}\bigg)\tilde{G}_{\frac{n-1}{2}+k,\frac{n-1}{2}}\bigg(\frac{nx}{s(1-\rho^2)}\bigg),\quad x\geq0, \\
F_{\overline{Z}_n}(x)&=\frac{(1-\rho^2)^{n/2}}{2\sqrt{\pi}\Gamma(n/2)}\sum_{k=0}^\infty\frac{(-2\rho)^k}{k!}\Gamma\bigg(\frac{k+1}{2}\bigg)\bigg(\frac{n+k}{2}\bigg)\tilde{G}_{\frac{n-1}{2}+k,\frac{n-1}{2}}\bigg(-\frac{nx}{s(1-\rho^2)}\bigg), \quad x<0,
\end{align*}
and, for $x\in\mathbb{R}$,
\begin{align*}F_{\overline{Z}_n}(x)&=\frac{1}{2}-\frac{\Gamma((n+1)/2)}{\sqrt{\pi}\Gamma(n/2)}\rho(1-\rho^2)^{n/2}{}_2F_1\bigg(1,\frac{n+1}{2};\frac{3}{2};\rho^2\bigg)\\
&\quad+\frac{(1-\rho^2)^{n/2}}{2\sqrt{\pi}\Gamma(n/2)}\sum_{k=0}^\infty(\mathrm{sgn}(x))^{k+1}\frac{(2\rho)^k}{k!}\Gamma\bigg(\frac{k+1}{2}\bigg)\bigg(\frac{n+k}{2}\bigg)G_{\frac{n-1}{2}+k,\frac{n-1}{2}}\bigg(\frac{n|x|}{s(1-\rho^2)}\bigg).
\end{align*}
In particular,
\begin{equation}\label{pz0}\mathbb{P}(\overline{Z}_n\leq 0)=\frac{1}{2}-\frac{\Gamma((n+1)/2)}{\sqrt{\pi}\Gamma(n/2)}\rho(1-\rho^2)^{n/2}{}_2F_1\bigg(1,\frac{n+1}{2};\frac{3}{2};\rho^2\bigg).
\end{equation}
\end{corollary}

\begin{remark}On setting $n=1$ in (\ref{pz0}) and using the formula (\ref{sin}), we obtain
\begin{equation*}\mathbb{P}(Z\leq0)=\frac{1}{2}-\frac{1}{\pi}\sin^{-1}(\rho),
\end{equation*}
which can also be deduced from the standard result that $\mathbb{P}(U\leq0, V>0)=\mathbb{P}(U>0, V\leq0)=1/4-\sin^{-1}(\rho)/(2\pi)$, for $(U,V)$ a bivariate normal random vector as defined above.
\end{remark}



\appendix

\section{Special functions}\label{appa}
In this appendix, we define the modified Bessel function of the second kind, the modified Lommel function of the first kind and the hypergeometric function, and present some basic properties that are used in this paper. Unless otherwise stated, the properties listed below can be found in \cite{olver}. For the modified Lommel function of the first kind, formulas (\ref{struve}), (\ref{tend0}) and (\ref{infty}) are given in \cite{gaunt lommel}, the integral formula (\ref{dfg}) can be found in \cite{r64}, whilst the results in (\ref{kint2})--(\ref{gineq}) are simple deductions from other properties listed in this appendix.

The \emph{modified Bessel function of the second kind} $K_\nu(x)$ is defined, for $\nu\in\mathbb{R}$ and $x>0$, by
\[K_\nu(x)=\int_0^\infty \mathrm{e}^{-x\cosh(t)}\cosh(\nu t)\,\mathrm{d}t.
\]
The \emph{generalized hypergeometric function} is defined by the power series
\begin{equation}
\label{gauss}
{}_pF_q(a_1,\ldots,a_p; b_1,\ldots,b_q;x)=\sum_{j=0}^\infty\frac{(a_1)_j\cdots(a_p)_j}{(b_1)_j\cdots(b_q)_j}\frac{x^j}{j!},
\end{equation}
for $|x|<1$, and by analytic continuation elsewhere. Here $(u)_j=u(u+1)\cdots(u+k-1)$ is the ascending factorial. The function ${}_2F_1(a,b;c;x)$ is known as the \emph{(Gaussian) hypergeometric function}. We have the special case 
\begin{equation}\label{sin} {}_2F_1(a,b;c;x)=\frac{\sin^{-1}(\sqrt{x})}{\sqrt{x(1-x)}}
\end{equation}
(see \texttt{http://functions.wolfram.com/07.23.03.3098.01}).

The \emph{modified Lommel function of the first kind} is defined by the hypergeometric series
\begin{align}t_{\mu,\nu}(x)&=\frac{x^{\mu+1}}{(\mu-\nu+1)(\mu+\nu+1)} {}_1F_2\bigg(1;\frac{\mu-\nu+3}{2},\frac{\mu+\nu+3}{2};\frac{x^2}{4}\bigg) \nonumber\\
&=2^{\mu-1}\Gamma\bigg(\frac{\mu-\nu+1}{2}\bigg)\Gamma\bigg(\frac{\mu+\nu+1}{2}\bigg)\sum_{k=0}^\infty\frac{(\frac{1}{2}x)^{\mu+2k+1}}{\Gamma\big(k+\frac{\mu-\nu+3}{2}\big)\Gamma\big(k+\frac{\mu+\nu+3}{2}\big)}. \nonumber
\end{align}
In this paper, it will be convenient to work with the following normalisation of the modified Lommel function of the first kind that was introduced by \cite{gaunt lommel}:
\begin{align*}\tilde{t}_{\mu,\nu}(x)&=\frac{1}{2^{\mu-1}\Gamma\big(\frac{\mu-\nu+1}{2}\big)\Gamma\big(\frac{\mu+\nu+1}{2}\big)}t_{\mu,\nu}(x) \\
&=\frac{1}{2^{\mu+1}\Gamma\big(\frac{\mu-\nu+3}{2}\big)\Gamma\big(\frac{\mu+\nu+3}{2}\big)} {}_1F_2\bigg(1;\frac{\mu-\nu+3}{2},\frac{\mu+\nu+3}{2};\frac{x^2}{4}\bigg).
\end{align*}

For $\nu=m+1/2$, $m=0,1,2,\ldots$, the modified Bessel function of the second kind takes an elementary form:
\begin{equation}\label{special} K_{m+1/2}(x)=\sqrt{\frac{\pi}{2x}}\sum_{j=0}^m\frac{(m+j)!}{(m-j)!j!}(2x)^{-j}\mathrm{e}^{-x}.
\end{equation}
The modified Struve function of the first kind $\mathbf{L}_\nu(x)$ is a special case of the function $\tilde{t}_{\mu,\nu}(x)$:
\begin{equation}\label{struve}\tilde{t}_{\nu,\nu}(x)=\mathbf{L}_\nu(x).
\end{equation}
The functions $K_\nu(x)$ and $\tilde{t}_{\mu,\nu}(x)$ have the following asymptotic behaviour:
\begin{eqnarray}\label{Ktend0}K_{\nu} (x) &\sim& \begin{cases} 2^{|\nu| -1} \Gamma (|\nu|) x^{-|\nu|}, & \quad x \downarrow 0, \: \nu \not= 0, \\
-\log x, & \quad x \downarrow 0, \: \nu = 0, \end{cases} \\
\label{Ktendinfinity} K_{\nu} (x) &\sim& \sqrt{\frac{\pi}{2x}} \mathrm{e}^{-x}, \quad x \rightarrow \infty,\: \nu\in\mathbb{R}, \\
\label{tend0}\tilde{t}_{\mu,\nu}(x)&\sim& \frac{(\frac{1}{2}x)^{\mu+1}}{\Gamma\big(\frac{\mu-\nu+3}{2}\big)\Gamma\big(\frac{\mu+\nu+3}{2}\big)}, \quad x\downarrow0,\:\mu>-3,\:|\nu|<\mu+3, \\
\label{infty}\tilde{t}_{\mu,\nu}(x)&\sim& \sim\frac{\mathrm{e}^x}{\sqrt{2\pi x}}, \quad x\rightarrow\infty, \:\mu,\nu\in\mathbb{R}.
\end{eqnarray}
The functions $K_\nu(x)$ and $\tilde{t}_{\mu,\nu}(x)$ are linked through the indefinite integral formula
\begin{align}\label{dfg}\int x^\mu K_\nu(x)\,\mathrm{d}x=-2^{\mu-1}\Gamma\bigg(\frac{\mu-\nu+1}{2}\bigg)\Gamma\bigg(\frac{\mu+\nu+1}{2}\bigg)G_{\mu,\nu}(x),
\end{align}
where $G_{\mu,\nu}(x)$ is defined as in (\ref{gmn}). With this indefinite integral formula and the limiting forms (\ref{Ktend0})--(\ref{infty}), we deduce the following integral formulas. For $\mu\geq\nu>-1/2$, $a>0$ and $x>0$,
\begin{align}\label{kint2}\int_0^x t^\mu K_\nu(at)\,\mathrm{d}t&=\frac{2^{\mu-1}}{a^\mu}\Gamma\bigg(\frac{\mu-\nu+1}{2}\bigg)\Gamma\bigg(\frac{\mu+\nu+1}{2}\bigg)G_{\mu,\nu}(ax), \\
\label{kint}\int_x^\infty t^\mu K_\nu(at)\,\mathrm{d}t&=\frac{2^{\mu-1}}{a^\mu}\Gamma\bigg(\frac{\mu-\nu+1}{2}\bigg)\Gamma\bigg(\frac{\mu+\nu+1}{2}\bigg)\tilde{G}_{\mu,\nu}(ax), 
\end{align}
where $\tilde{G}_{\mu,\nu}(x)$ is defined as in (\ref{gmn2}).

 Since $K_\nu(x)>0$ for all $\nu\in\mathbb{R}$, $x>0$, and the gamma functions in (\ref{kint2}) and (\ref{kint}) are positive for $\mu\geq\nu>-1/2$, it follows that, for fixed $\mu\geq\nu>-1/2$,   $G_{\mu,\nu}(x)$ is an increasing function of $x$ on $(0,\infty)$ with $G_{\mu,\nu}(x)>0$, and $\tilde{G}_{\mu,\nu}(x)$ is a decreasing function of $x$ on $(0,\infty)$ with $\tilde{G}_{\mu,\nu}(x)>0$. Therefore, since $\tilde{G}_{\mu,\nu}(x)=1-G_{\mu,\nu}(x)$, we deduce that, for  $\mu\geq\nu>-1/2$, $x>0$,
\begin{align}\label{gineq}0<G_{\mu,\nu}(x)<1, \quad 0<\tilde{G}_{\mu,\nu}(x)<1.
\end{align}


\end{document}